\documentclass{article}
\usepackage{color}
\catcode`\@=11 \@addtoreset{equation}{section}

\catcode`\@=12

\newtheorem{Theorem}{Theorem}[section]
\newtheorem{Proposition}{Proposition}[section]
\newtheorem{Lemma}{Lemma}[section]
\newtheorem{Corollary}{Corollary}[section]
\newtheorem{Remark}{Remark}[section]

\newcommand{\bTheorem}[1]{
\begin{Theorem} \label{T#1} }
\newcommand{\eT}{\end{Theorem}}

\newcommand{\bProposition}[1]{
\begin{Proposition} \label{P#1}}
\newcommand{\eP}{\end{Proposition}}

\newcommand{\bLemma}[1]{
\begin{Lemma} \label{L#1} }
\newcommand{\eL}{\end{Lemma}}

\newcommand{\bCorollary}[1]{
\begin{Corollary} \label{C#1} }
\newcommand{\eC}{\end{Corollary}}

\newcommand{\bFormula}[1]{
\begin{equation} \label{#1}}
\newcommand{\eF}{\end{equation}}

\newcommand{\bRemark}[1]{
\begin{Remark} \label{R#1} }
\newcommand{\eR}{\end{Remark}}

\newcommand{\Ov}[1]{\overline{#1}}
\newcommand{\DC}{C^\infty_c}

\newcommand{\vr}{\varrho}
\newcommand{\vre}{\vr_\ep}
\newcommand{\vte}{\vt_\ep}
\newcommand{\vue}{\vu_\ep}

\newcommand{\vt}{\vartheta}
\newcommand{\vu}{\vc{u}}
\newcommand{\vc}[1]{{\bf #1}}

\newcommand{\Div}{{\rm div}_x}
\newcommand{\Grad}{\nabla_x}

\newcommand{\tn}[1]{\mbox {\F #1}}
\newcommand{\dx}{{\rm d} {x}}
\newcommand{\dt}{{\rm d} t }

\newcommand{\dxdt}{\dx \ \dt}
\newcommand{\intO}[1]{\int_{\Omega} #1 \ \dx}

\newcommand{\ep}{\varepsilon}

\font\F=msbm10 scaled 1000

\definecolor{grey}{rgb}{0.85,0.85,0.85}
\date{}

\long\def\greybox#1{%
    \newbox\contentbox%
    \newbox\bkgdbox%
    \setbox\contentbox\hbox to \hsize{%
        \vtop{
            \kern\columnsep
            \hbox to \hsize{%
                \kern\columnsep%
                \advance\hsize by -2\columnsep%
                \setlength{\textwidth}{\hsize}%
                \vbox{
                    \parskip=\baselineskip
                    \parindent=0bp
                    #1
                }%
                \kern\columnsep%
            }%
            \kern\columnsep%
        }%
    }%
    \setbox\bkgdbox\vbox{
        \color{grey}
        \hrule width  \wd\contentbox %
               height \ht\contentbox %
               depth  \dp\contentbox
        \color{black}
    }%
    \wd\bkgdbox=0bp%
    \vbox{\hbox to \hsize{\box\bkgdbox\box\contentbox}}%
    \vskip\baselineskip%
}

\begin{document}


\title{Inviscid incompressible limits under mild stratification: A rigorous derivation of the Euler-Boussinesq system}
\author{Eduard Feireisl \thanks{The research of E.F. leading to these results has received funding from the European Research Council under the European Union's Seventh Framework
Programme (FP7/2007-2013)/ ERC Grant Agreement 320078.} \and Anton\' \i n Novotn\' y \thanks{The work was supported by the MODTERCOM project within the APEX programme of the region Provence-Alpe-C\^ote d'Azur and by  RVO: 67985840}}

\maketitle

\bigskip

\centerline{Institute of Mathematics of the Academy of Sciences of
the Czech Republic}

\centerline{\v Zitn\' a 25, 115 67 Praha 1,
Czech Republic}

\bigskip

\centerline{IMATH, Universit\' e du Sud Toulon-Var}

\centerline{BP 20139, 839 57 La Garde, France}

\medskip

\begin{abstract}
We consider the full Navier-Stokes-Fourier system in the singular regime of small
Mach and large Reynolds and P{\' e}clet numbers, with ill prepared initial data on
an unbounded domain$\Omega \subset R^3$ with a compact boundary. We perform the singular limit in the framework of
weak solutions and identify
the Euler-Boussinesq system as the target problem.
\end{abstract}

\tableofcontents

\section{Introduction}
\label{i}

The present paper is an extension of our previous results concerning the inviscid incompressible limit
of the Navier-Stokes-Fourier system \cite{FeiNov12}. In contrast with \cite{FeiNov12}, where the problem is considered on the whole space $R^3$
without any driving force imposed, we consider a more realistic situation when the fluid is subject to a gravitational force due to the physical objects placed
\emph{outside} the fluid domain. Accordingly, we shall assume that the fluid occupies an unbounded \emph{exterior} domain $\Omega \subset R^3$ with smooth
(compact) boundary. Such a situation is interesting from the point of view of possible applications in various meteorological models as the singular limit in the
low Mach, Froude, and large Reynolds and P\' eclet numbers leads to a target system driven by the buoyancy force proportional to temperature deviations.
In particular, we provide a rigorous justification of the so-called Euler-Boussinesq approximation. Our approach is based on the recently discovered
relative entropy inequality \cite{FeiNov10} and the related concept of \emph{dissipative solution} for the Navier-Stokes-Fourier system.
In comparison with \cite{FeiNov12}, the present problem features some additional mathematical difficulties related to the geometry of the underlying spatial domain
and the presence of a driving force. In particular, we have to handle perturbations of weakly
stratified equilibrium states, whereas those are simply constant in \cite{FeiNov12}.

We consider the motion of a compressible, viscous and heat conducting fluid, with the density $\vr = \vr(t,x)$, the velocity $\vu = \vu(t,x)$, and
the absolute temperature
$\vt = \vt(t,x)$ governed by the scaled \emph{Navier-Stokes-Fourier system}:

\bFormula{i1}
\partial_t \vr + \Div (\vr \vu) = 0,
\eF
\bFormula{i2}
\partial_t (\vr \vu) + \Div (\vr \vu \otimes \vu) + \frac{1}{\ep^2} \Grad p(\vr, \vt) = \ep^a \Div \tn{S} (\vt, \Grad \vu) +{\frac 1\ep}\vr \Grad F,
\eF
\bFormula{i3}
\partial_t (\vr s(\vr , \vt)) + \Div (\vr s(\vr, \vt) \vu) + \ep^\beta \Div \left( \frac{\vc{q}(\vt, \Grad \vt)}{\vt} \right) =
\frac{1}{\vt} \left( \ep^{2 + a} \tn{S} (\vt, \Grad \vu) : \Grad \vu - \ep^b \frac{\vc{q}(\vt, \Grad \vt) \cdot
\Grad \vt }{\vt} \right),
\eF
where $p = p(\vr,\vt)$ is the pressure, $s = s(\vr, \vt)$ the specific entropy, the symbol $\tn{S}(\vt, \Grad \vu)$ denotes
the viscous stress satisfying \emph{Newton's law}
\bFormula{i3a}
\tn{S}(\vt, \Grad \vu) = \mu(\vt) \left( \Grad \vu + \Grad^t \vu - \frac{2}{3} \Div \vu \right) + \eta(\vt) \Div \vu \tn{I},
\eF
and $\vc{q} = \vc{q}(\vt, \Grad \vt)$ is the heat flux determined by \emph{Fourier's law}
\bFormula{i3b}
\vc{q}(\vt, \Grad \vt) = - \kappa (\vt) \Grad \vt,
\eF
where the quantities $\mu$, $\eta$, $\kappa$ are temperature dependent transport coefficients.

The fluid occupies an exterior domain $\Omega \subset R^3$, with impermeable, thermally insulating and frictionless boundary, specifically,
\bFormula{i3c}
\vu \cdot \vc{n} = [\tn{S}(\vt, \Grad \vu) \cdot \vc{n}]_{\rm tan}|_{\partial \Omega} = 0, \ \Grad \vt \cdot \vc{n}|_{\partial \Omega} = 0.
\eF
In addition, we consider the far field boundary conditions
\bFormula{i5}
\vr \to \Ov{\vr}, \ \vt \to \Ov{\vt}, \ \vu \to 0 \ \mbox{as} \ |x| \to \infty,
\eF
where $\Ov{\vr}$, $\Ov{\vt}$ are positive constants.

The scaling in (\ref{i1} - \ref{i3}), expressed by means of a single (small) parameter $\ep$, corresponds to:

\medskip

\noindent
Mach number \dotfill $\ep$,

\noindent
Froude number \dotfill $\ep^{1/2}$,

\noindent
Reynolds number \dotfill $\ep^{-a}$,

\noindent
P\' eclet number \dotfill $\ep^{-b}$.

\medskip

In accordance with the previous discussion, we consider the driving force induced by a potential
\bFormula{potential1}
F(x) = \int_{R^3} \frac{1}{x - y} m(y) {\rm d}y,\ m \geq 0, \ {\rm supp}[m] \subset R^3 \setminus \Omega,
\eF
meaning the fluid is driven by the gravitational force of objects lying outside the fluid domain.

Finally, the initial data are taken in the form
\bFormula{i4}
\vr(0, \cdot) = \vr_{0,\ep} = \Ov{\vr}_\ep + \ep \vr^{(1)}_{0,\ep},\
\vt(0, \cdot) = \vt_{0,\ep} = \Ov{\vt} + \ep \vt^{(1)}_{0,\ep},\ \vu(0, \cdot) = \vu_{0, \ep},
\eF
where $(\overline\vr_\ep, \overline\vt)$ is the equilibrium solution associated with the far field values of $\Ov{\vr}$, $\Ov{\vt}$, namely
\begin{equation}\label{equilibrium}
\Grad p(\overline\vre,\overline\vt) = \ep\overline\vre\Grad F,\quad\overline\vre\to\overline\vr\quad\mbox{as $|x|\to \infty$}.
\end{equation}

The limit (target) problem can be formally identified as the incompressible \emph{Euler-Boussinesq system}:
\bFormula{l1}
\Div \vc{v} = 0,
\eF
\bFormula{l2}
\partial_t \vc{v} + \vc{v} \cdot \Grad \vc{v} + \Grad \Pi  = - a (\Ov{\vr}, \Ov{\vt} )\theta\Grad F,
\eF
\bFormula{l3}
c_p(\Ov{\vr}, \Ov{\vt}) \left( \partial_t \theta + \vc{v} \cdot \Grad \theta \right) - \Ov{\vt} a(\Ov{\vr}, \Ov{\vt}) \vc{v} \cdot \Grad F = 0,
\eF
where we have denoted

\medskip

\noindent
thermal expansion coefficient \dotfill $a(\Ov{\vr}, \Ov{\vt})$,

\noindent
specific heat at constant pressure \dotfill $c_p(\Ov{\vr}, \Ov{\vt})$,

\medskip

\noindent
cf.  \cite[Chapter 5]{FeNo6} and \cite{FeiNov12}.
Here, the function $\vc{v}$ is the limit velocity, while $\theta$ is associated with the asymptotic temperature
(entropy) deviations
\[
\theta \approx \frac{\vte - \Ov{\vt}}{\ep}.
\]
The
exact statement of our results including the initial data for the target system (\ref{l1} - \ref{l3}) will be specified in Theorem \ref{Tm1} below.

We address the problem in the framework of weak solutions for the Navier-Stokes-Fourier system (\ref{i1} - \ref{i3}),
developed in \cite{FeNo6}, and later extended to problems on unbounded domains in \cite{JeJiNo}. The main advantage of this approach is
the convergence towards the target system on any time interval $[0,T]$, on which the Euler-Boussinsesq system (\ref{l1}), (\ref{l2})
possesses a regular solution. We refer to Masmoudi \cite{MAS5} for related results on the compressible barotropic Navier-Stokes system
in the whole space $R^3$, see also the survey \cite{MAS1}. An alternative approach to singular limits, proposed in the seminal paper by
Klainerman and Majda \cite{KM1}, uses the strong solutions for both the primitive and the target system that may exist, however, only
on a possible very short time interval. Using the { same} framework,
Alazard \cite{AL2}, \cite{AL1}, \cite{AL} addresses several
singular limits of the compressible Euler and/or Navier-Stokes-Fourier system, in the absence of external forcing. The present setting, where
the action of the gravitation gives rise to the buoyancy force proportional to $-\theta\Grad F$, represents a stronger coupling between the equations,
typical for certain models used in meteorology and physics of the atmosphere, see Klein \cite{KL1}, \cite{Klein}, Zeytounian \cite{ZEY1}.

The necessary preliminary material including various concepts of weak solutions to
the Navier-Stokes-Fourier system is collected in  Section \ref{v}. Section \ref{m} contains the main result on the asymptotic limit for $\ep \to 0$, the
proof of which is the main objective of the remaining part for the paper. In Section \ref{b}, the relative entropy inequality is
 used to establish the necessary uniform bounds independent of $\ep \to 0$. The problem of propagation and dispersion of the associated
acoustic waves is discussed in Section \ref{r}. The proof of convergence towards the limit system is completed in Section \ref{c}.

\section{Preliminaries, weak solutions to the Navier-Stokes-Fourier system}
\label{v}

Motivated by \cite{FeiNov10}, we introduce the \emph{relative entropy functional}
\bFormula{r0}
\mathcal{E}_\ep \left( \vr, \vt, \vu \Big| r , \Theta, \vc{U} \right) =
\intO{ \left[ \frac{1}{2} \vr |\vu - \vc{U} |^2 + \frac{1}{\ep^2} \left( H_\Theta (\vr, \vt) -
\frac{\partial H_\Theta (r, \Theta)}{\partial\vr}(\vr - r)- H_\Theta (r, \Theta) \right)
\right] },
\eF
where
\bFormula{bfe}
H_\Theta (\vr, \vt) = \vr \Big( e(\vr, \vt) - \Theta s(\vr,\vt) \Big)
\eF
is the ballistic free energy. We say that a trio of functions $\{ \vr, \vt, \vu \}$ represents a
\emph{dissipative weak solution} of the Navier-Stokes-Fourier system (\ref{i1} - \ref{i5}) in $(0,T) \times \Omega$ if:

\begin{itemize}

\item $\vr \geq 0$, $\vt > 0$ a.a. in $(0,T) \times \Omega$,
\[
( \vr - \Ov{\vr}_\ep ) \in L^\infty(0,T; L^2 + L^{5/3}(\Omega)),\ (\vt - \Ov{\vt}) \in L^\infty(0, T; L^2 + L^4 (\Omega)),
\]
\[
\Grad \vt , \ \Grad \log(\vt) \in L^2(0,T; L^2(\Omega;R^3)),
\]
\[
\vu \in L^2(0,T; W^{1,2}(\Omega; R^3)),\ \vu \cdot \vc{n}|_{\partial \Omega} = 0,
\]
where $[\Ov{\vr}_\ep, \Ov{\vt}]$ stands for the equilibrium solution introduced in (\ref{equilibrium});

\item the equation of continuity (\ref{i1}) holds as a family of integral identities
\bFormula{v1}
\int_{\Omega} \Big[ \vr(\tau, \cdot) \varphi (\tau, \cdot) - \vr_{0,\ep} \varphi (0, \cdot) \Big] \ \dx
= \int_0^\tau \int_{R^3} \Big( \vr \partial_t \varphi + \vr \vu \cdot \Grad \varphi \Big) \ \dxdt
\eF
for any $\tau \in [0,T]$ and any test function $\varphi \in \DC([0,T] \times \Ov{\Omega})$;

\item the momentum equation (\ref{i2}), together with the initial condition (\ref{i4}), are satisfied in the sense of distributions,
\bFormula{v2}
\int_{\Omega} \Big[ \vr \vu (\tau, \cdot) \cdot \varphi (\tau, \cdot) - \vr_{0,\ep} \vu_{0,\ep}
\varphi(0, \cdot) \Big] \ \dx
\eF
\[
= \int_0^\tau \int_{\Omega} \Big( \vr \vu \cdot \partial_t \varphi + \vr \vu \otimes \vu : \Grad \varphi +
\frac{1}{\ep^2} p(\vr,\vt) \Div \varphi - \ep^{a} \tn{S}(\vt, \Grad \vu) : \Grad \varphi + \frac{1}{\ep} \Grad F \cdot \varphi \Big) \ \dxdt
\]
for any $\tau \in [0,T]$, and any $\varphi \in \DC([0,T] \times \Ov{\Omega}; R^3)$, $\varphi \cdot \vc{n}|_{\partial \Omega} = 0$;

\item the entropy production equation (\ref{i3}) is relaxed to the entropy inequality
\bFormula{v3}
\int_{\Omega} \Big[ \vr_{0,\ep} s(\vr_{0,\ep}, \vt_{0, \ep} ) \varphi(0, \cdot) -
\vr s(\vr, \vt) (\tau, \cdot) \varphi(\tau, \cdot) \Big] \ \dx
\eF
\[
+ \int_0^\tau \int_{\Omega} \frac{1}{\vt} \left( \ep^{2 + a} \tn{S}(\vt, \Grad \vu) : \Grad \vu - \ep^b
\frac{\vc{q}(\vt, \Grad \vt) \cdot \Grad \vt }{\vt} \right) \varphi \ \dxdt
\]
\[
\leq - \int_0^\tau \int_{\Omega} \left( \vr s(\vr,\vt) \partial_t \varphi + \vr s(\vr, \vt) \vu \cdot \Grad \varphi +
\ep^b \frac{ \vc{q}(\vt, \Grad \vt) }{\vt} \cdot \Grad \varphi \right) \ \dxdt
\]
for a.a. $\tau \in [0,T]$ and any test function $\varphi \in \DC([0,T] \times \Ov{\Omega})$, $\varphi \geq 0$;

\item the \emph{relative entropy inequality}

\bFormula{r1}
\left[ \mathcal{E}_\ep \left( \vr, \vt, \vu \Big| r , \Theta, \vc{U} \right) \right]_{t = 0}^\tau
+ \int_0^\tau \intO{ \frac{\Theta}{\vt} \left( \ep^a \tn{S} (\vt, \Grad \vu) : \Grad \vu -
\ep^{b-2} \frac{\vc{q}(\vt, \Grad \vt) \cdot \Grad \vt }{\vt} \right) } \ \dt
\eF
\[
\leq \int_0^\tau \intO{ \Big( \vr \Big( \partial_t \vc{U} +
\vu\cdot \Grad \vc{U} \Big) \cdot (\vc{U} - \vu)  + \ep^a
\tn{S}(\vt, \Grad \vu): \Grad \vc{U} \Big) } \ \dt
\]
\[
+\frac
1{\ep^2}\int_0^\tau\intO{\Big[\Big(p(r,\Theta)-p(\vr,\vt)\Big){\rm
div}\vc U +\frac\vr {r}(\vc U-\vu)\cdot\Grad
p(r,\Theta)\Big]}{\rm d}t
\]
\[
- \frac{1}{\ep^2} \int_0^\tau \intO{ \left( \vr \Big( s(\vr,\vt) - s(r, \Theta) \Big) \partial_t \Theta +
\vr \Big( s(\vr,\vt) - s(r, \Theta) \Big) \vu \cdot \Grad \Theta + \ep^b \frac{\vc{q}(\vt, \Grad \vt) }{\vt} \cdot
\Grad \Theta \right) } \ \dt
\]
\[
+ \frac{1}{\ep^2} \int_0^\tau \intO{ \frac{r - \vr}
{r}\Big(
\partial_t p(r, \Theta) + \vc U \cdot \Grad p(r,
\Theta)\Big)  } \ \dt { -\frac 1\ep\int_0^\tau\int_{\Omega}\vre\Grad F\cdot (\vc U_\ep-\vue){\rm d}x}.
\]
holds for a.a. $t\in (0,T)$ and for
any trio of continuously differentiable ``test'' functions defined on $[0,T] \times \Ov{\Omega}$,
\[
r > 0 , \ \Theta > 0 , \ r \equiv \Ov{\vr} , \ \Theta \equiv \Ov{\vt} \ \mbox{outside a compact subset of}\ \Ov{\Omega},
\]
\[
\vc{U} \in C([ 0, T]; W^{k,2}(\Omega;R^3)),\ \partial_t \vc{U} \in \ C([0,T]; W^{k-1,2}(\Omega;R^3)),\ k > \frac{5}{2}, \ \vc{U} \cdot \vc{n}|_{\partial
\Omega} = 0.
\]

\end{itemize}

\bRemark{i1}
Note that the above definition of dissipative weak solutions on \emph{unbounded} domains, proposed in \cite{JeJiNo}, is different from
that on bounded domains introduced in \cite{FeiNov10}. In \cite{FeiNov10}, the relative entropy inequality (\ref{r1}) is replaced by the total
energy balance, whereas (\ref{r1}) is automatically satisfied for any weak solution to the Navier-Stokes-Fourier system. The
weak solutions introduced in this paper can be therefore viewed as ``very weak dissipative solutions'' of the primitive system.
\eR

\subsection{Structural restrictions imposed on constitutive relations}

We study our singular limit problem under certain physically motivated restrictions imposed on constitutive equations. They are basically the same
as required by the existence theory developed in \cite[Chapter 3]{FeNo6}. Although they might be slightly relaxed if only the convergence towards the target
system is studied, we list them in the form presented in  \cite[Chapter 3]{FeNo6}, where the interested reader may find more information
concerning the physical background as well as possible generalizations.

The pressure $p=p(\vr,\vt)$ is given by the formula
\bFormula{i10}
p(\vr, \vt) = \vt^{5/2} P \left( \frac{\vr}{\vt^{3/2}} \right) + \frac{a}{3} \vt^4 , \ a  > 0;
\eF
the specific internal energy $e = e(\vr,\vt)$ and the specific entropy $s = s(\vr, \vt)$ read
\bFormula{i11}
e(\vr, \vt) = \frac{3}{2} \vt \frac{ \vt^{3/2} }{\vr} P \left( \frac{\vr}{\vt^{3/2}} \right) + {a} \vt^4
\eF
\bFormula{i12}
s(\vr, \vt) = S \left( \frac{\vr}{\vt^{3/2}} \right) + \frac{4a}{3} \frac{\vt^3}{\vr},
\eF
where
\bFormula{i13}
P \in C^1[0, \infty) \cap C^3(0,\infty), \ P(0) = 0 ,\ P'(Z) > 0 \ \mbox{for all}\ Z \geq 0,
\eF
\bFormula{i14}
\lim_{Z \to \infty} \frac{P(Z)}{Z^{5/3}} = P_\infty > 0,
\eF
\bFormula{i15}
0 < \frac{ \frac{5}{3} P(Z) - P'(Z) Z }{Z} < c \ \mbox{for all} \ Z > 0,
\eF
and
\bFormula{i16}
S'(Z) = - \frac{3}{2} \frac{ \frac{5}{3} P(Z) - P'(Z) Z }{Z^2} , \ \lim_{Z \to \infty} S(Z) = 0.
\eF
The relation (\ref{i15}) expresses  positivity and uniform boundedness of the specific heat at constant volume.

The transport coefficients $\mu$, $\eta$, and $\kappa$ are effective functions of the temperature,
\bFormula{i8}
\mu, \ \eta  \in C^1[0, \infty) \ \mbox{are globally Lipschitz in,} \ [0, \infty), \ 0 < \underline{\mu} (1 + \vt) \leq \mu(\vt),\
\eta(\vt) \geq 0,
\ \mbox{for all}\ \vt \geq 0,
\eF
\bFormula{i9}
\kappa \in C^1[0, \infty),\ 0 < \underline{\kappa}(1 + \vt^3) \leq \kappa (\vt) \leq \Ov{\kappa}(1 + \vt^3) \ \mbox{for all}\ \vt \geq 0.
\eF

\subsection{Target system}

As noted in the introduction, the expected limit is the Euler-Boussinesq system (\ref{l1}- \ref{l3}) endowed with  the
initial data
\begin{equation}\label{inittarget}
\theta_0(0,\cdot)=\theta_0,\;\vc{v}(0, \cdot) = \vc{v}_0.
\end{equation}

In agreement with the nowadays standard theory of well-posedness for hyperbolic systems, see e.g. Kato \cite{Kato},
we suppose that the system (\ref{l1}- \ref{l3}), endowed with the initial data
\begin{equation}\label{initdata}
(\theta_0, \vc{v}_0)\in  W^{k,2}(\Omega;R^4), \ \|(\theta_0,\vc{v}_0)\|_{W^{k,2}(\Omega;R^4)}\le D, \ \Div \vc{v}_0 = 0,\
\vc{v}_0 \cdot \vc{n}|_{\partial \Omega} = 0, \ k > \frac{5}{2},
\end{equation}
possesses a regular solution $(\theta,\vc{v})$,
\bFormula{euler} (\theta,\vc{v}) \in C([0, T_{\rm max}); W^{k,2}(\Omega;R^4)),
\ (\partial_t \vc{v},\,\Grad \Pi)\in C([0, T_{\rm max}); W^{k-1,2}(\Omega;R^6)),
\eF defined on a maximal time interval $[0, T_{\rm max})$, $T_{\rm
max} = T_{\rm max}(D)$.

\subsection{Equilibrium state}

We finish this preliminary part by recalling the basic properties of the equilibrium solution $(\Ov{\vr}_\ep, \Ov{\vt})$.
Since the potential $F$ is given by (\ref{potential1}), it is easy to check that
\bFormula{p2} \partial_\vr H_{\overline\vt}(\overline\vre,\overline\vt) = \ep F +
\partial_\vr H_{\overline\vt}(\overline\vr,\overline\vt); \eF
whence, under the assumptions (\ref{i10}), (\ref{i13}--\ref{i15}),
\bFormula{p3}
\overline\vre\in C^3(\Omega),\
\left| \frac{\overline{\vr}_\ep (x) - \overline{\vr}}{\ep} \right| \leq c  F (x) \quad \left| \Grad \vre(x)  \right| \leq \ep |\Grad F(x)| ,\ x \in \Omega.
\eF
The reader can consult \cite{FeiSch} for details.

\section{Main result}

\label{m}

For a vector field $\vc{U} \in L^2(\Omega;R^3)$, we denote by $\vc{H}[\vc{U}]$ the standard \emph{Helmholtz projection}
on the space of solenoidal functions.

We are ready to state the main result of this paper.

\bTheorem{m1}
Let the thermodynamic functions
$p$, $e$, $s$, and the transport coefficients $\mu$, $\eta$, $\kappa$ satisfy the hypotheses (\ref{i10} - \ref{i16}),
(\ref{i8}), (\ref{i9}). Let the potential force $F$ be given by (\ref{potential1}). Let the exponents $a,b$, determining the
Reynold and P\'eclet number scales, satisfy
\bFormula{cof}
b > 0, \ 0 < a < \frac{10}{3}.
\eF

Next, let
the initial data (\ref{i4}) be chosen in such a way that
\bFormula{m1}
\{ \vr^{(1)}_{0, \ep} \}_{\ep > 0},\  \{ \vt^{(1)}_{0,\ep} \}_{\ep > 0} \ \mbox{are bounded in}\
L^2 \cap L^\infty (\Omega), \ \vr^{(1)}_{0, \ep} \to \vr^{(1)}_0,\ \vt^{(1)}_{0, \ep} \to \vt^{(1)}_0
\ \mbox{in}\ L^2(\Omega),
\eF
\bFormula{m2}
\{ \vu_{0,\ep} \}_{\ep > 0} \ \mbox{is bounded in}\ L^2(\Omega;R^3), \ \vu_{0, \ep} \to \vu_0 \ \mbox{in}\ L^2(\Omega;R^3),
\eF
where
\bFormula{HYP}
\vr^{(1)}_0 , \ \vt^{(1)}_0 \in W^{1,2} \cap W^{1,\infty}(\Omega),\ \vc{H}[\vu_0] = \vc{v}_0 \in W^{k,2}(\Omega;R^3) \ \mbox{for a certain}\ k > \frac{5}{2}.
\eF

Suppose that the Euler-Boussinesq system (\ref{l1}--\ref{l3}), endowed with the initial data
\bFormula{m6}
\vc{v}_0 = \vc{H}[ \vu_0 ], \ \theta_0 = \frac{\Ov{\vt}}{c_p(\Ov{\vr}, \Ov{\vt})} \left( \frac{\partial s(\Ov{\vr}, \Ov{\vt})}{
\partial \vr} \vr^{(1)}_0 + \frac{\partial s(\Ov{\vr}, \Ov{\vt})}{
\partial \vt} \vt^{(1)}_0 \right) ,
\eF
admits a regular solution $[\vc{v}, \theta]$ in the class (\ref{euler}) defined on a maximal time interval $[0, T_{\rm max})$.

Finally, let $\{ \vre, \vte, \vue \}$ be a dissipative weak solution of the Navier-Stokes-Fourier system (\ref{i1} - \ref{i5}) in
$(0,T) \times R^3$, $T < T_{\rm max}$.

Then
\bFormula{m3}
{\rm ess} \sup_{t \in (0,T)} \|  \vre (t, \cdot) - \Ov{\vr}  \|_{L^{5/3}_{\rm loc}(\Ov{\Omega})} \leq \ep c,
\eF
\bFormula{m4}
\sqrt{ {\vre} } \vue \to \sqrt{ \Ov{\vr} } \ \vc{v} \ \mbox{in}\ L^\infty_{\rm loc}((0,T]; L^2_{\rm loc} (\Ov{\Omega};R^3))
\ \mbox{and weakly-(*) in} \ L^\infty(0,T; L^2(\Omega;R^3)),
\eF
and
\bFormula{m5}
\frac{ \vte - \Ov{\vt} }{\ep} \to \theta \;
\mbox{in}\ L^\infty_{\rm loc}((0,T]; L^2_{\rm loc} (\Ov{\Omega})),
\ \mbox{and weakly-(*) in}\ L^\infty(0,T; L^2(\Omega)).
\eF

\eT
\bigskip\noindent

\bRemark{m1}
Under the hypotheses (\ref{i10} - \ref{i9}), the existence of
dissipative weak solutions to the Navier-Stokes-Fourier system in $(0,T)
\times \Omega$ was shown in \cite{JeJiNo}.
\eR

The rest of the paper is devoted to the proof of Theorem \ref{Tm1}.

\section{Uniform bounds}
\label{b}

In this section, we derive uniform bounds on the family of solutions $[\vre, \vue, \vte]$ \emph{independent} of the scaling parameter $\ep \to 0$.

\subsection{Energy bounds}

Taking  $r=\overline\vre$, $\Theta=\overline\vartheta$,
$\vc U=0$ as test functions in the relative entropy inequality (\ref{r1}) we obtain

\bFormula{v4}
\int_{\Omega} \left[ \frac{1}{2} \vre |\vue|^2 + \frac{1}{\ep^2} \left( H_{\Ov{\vt}} (\vre, \vte) -
\frac{\partial H_{\Ov{\vt}} (\Ov{\vr}_\ep, \Ov{\vt})}{\partial\vr}(\vre - \Ov{\vr}_\ep )- H_{\Ov{\vt}} (\Ov{\vr}_\ep, \Ov{\vt}) \right) \right]
{\rm d} x
\eF
\[
+ \Ov{\vt} \int_0^\tau \int_{\Omega} \frac{1}{\vte} \left( \ep^{a} \tn{S}(\vte, \Grad \vue) : \Grad \vue - \ep^{b - 2} \frac{\vc{q}(\vte, \Grad \vte) \cdot \Grad \vte
 }{\vt} \right) \ \dxdt
\]
\[
\leq \int_{\Omega} \left[ \frac{1}{2} \vr_{0,\ep} |\vu_{0,\ep}|^2 + \frac{1}{\ep^2} \left( H_{\Ov{\vt}} (\vr_{0,\ep}, \vt_{0,\ep}) -
\frac{\partial H_{\Ov{\vt}} (\Ov{\vr}_{\ep}, \Ov{\vt})}{\partial\vr}(\vr_{0,\ep} - \Ov{\vr}_\ep )- H_{\Ov{\vt}} (\Ov{\vr}_\ep, \Ov{\vt}) \right)\right]  \ \dx
\]
for a.a. $\tau \in [0,T]$. Note that such a choice of test functions can be justified by means of a density argument.
Thanks to the hypotheses (\ref{m1}), (\ref{m2}), the integral on the right-hand side of (\ref{v4}) remains bounded uniformly for $\ep \to 0$.

In accordance with the structural properties of the thermodynamic functions imposed through (\ref{i10} - \ref{i16}),
the ballistic free energy enjoys the following properties:
For any compact $K \subset (0, \infty)^2$ and $
(r, \Theta) \in K$,
there exists a strictly positive constant $c(K)$, depending only on $K$ and the structural properties of $P$, such that
\bFormula{b1}
\left( H_\Theta (\vr, \vt) -
\frac{\partial H_\Theta (r, \Theta)}{\partial\vr}(\vr - r)- H_\Theta (r, \Theta) \right) \geq c(K) \left( |\vr - r |^2 + |\vt - \Theta |^2 \right) \ \mbox{if}
\ (\vr, \vt) \in K,
\eF
\bFormula{b2}
\left( H_\Theta (\vr, \vt) -
\frac{\partial H_\Theta (r, \Theta)}{\partial\vr}(\vr - r)- H_\Theta (r, \Theta) \right) \ge c(K) \Big(1+\vr^\gamma +\vt^4\Big) \ \mbox{if}
\ (\vr, \vt) \in (0,\infty)^2 \setminus K.
\eF

Similarly to \cite[Chapter 4.7]{FeNo6}, we introduce a decomposition of a function $h$:
\[
h = [h]_{\rm ess} + [h]_{\rm res} \ \mbox{for a measurable function} \ h,
\]
where
\[
[h]_{\rm ess} = h \ 1_{ \{ \Ov{\vr} / 2 < \vre < 2 \Ov{\vr} ; \ \Ov{\vt}/2 < \vte < 2 \Ov{\vt} \} },\
[h]_{\rm res} = h - h_{\rm ess}.
\]

Combining (\ref{v4}), (\ref{b1}), (\ref{b2}), (\ref{p3}) with the hypotheses (\ref{i10} - \ref{i9}) we deduce the following estimates:
\bFormula{b3}
{\rm ess} \sup_{t \in (0,T)} \| \sqrt{\vre} \vue (t, \cdot) \|_{L^2(\Omega;R^3)} \leq c,
\eF
\bFormula{b4}
{\rm ess} \sup_{t \in (0,T)} \left\| \left[ \frac{\vre - \Ov{\vr}_\ep }{\ep}(t, \cdot) \right]_{\rm ess} \right\|_{L^2(\Omega;R^3)} +
{\rm ess} \sup_{t \in (0,T)} \left\| \left[ \frac{\vte - \Ov{\vt} }{\ep} (t, \cdot) \right]_{\rm ess} \right\|_{L^2(\Omega;R^3)} \leq c ,
\eF
\bFormula{b5}
{\rm ess} \sup_{t \in (0,T)} \int_{\Omega} \left( \left[ \vre^{5/3}(t, \cdot) \right]^{5/3}_{\rm res} + \left[ \vte (t, \cdot) \right]^4_{\rm res} + 1_{\rm res}(t, \cdot) \right) \ \dx \leq \ep^2 c,
\eF
and
\bFormula{b6}
\left\| \ep^{a/2} \vue \right\|_{L^2(0,T; W^{1,2}(\Omega;R^3))} \leq c,
\eF
\bFormula{b7}
\left\| \ep^{(b-2)/2} \left( \vte - \Ov{\vt} \right) \right\|_{L^2(0,T; W^{1,2}(\Omega;R^3))}
+ \left\| \ep^{(b-2)/2} \left( \log(\vte) - \log(\Ov{\vt}) \right) \right\|_{L^2(0,T; W^{1,2}(\Omega;R^3))} \leq c,
\eF
where the symbol $c$ stands for a generic constant independent of $\ep$. We remark that (\ref{b6}) follows from the generalized Korn's inequality
$
{ (\int_\Omega \vr_\ep\vc w^2{\rm d} x)^{1/2}+}\| \Grad \vc{w} + \Grad^t \vc{w} - \frac{2}{3} \Div \vc{w} \tn{I} \|_{L^2} \geq c
\| \Grad \vc{w}  \|_{L^2}$ for $\vc{w} \in W^{1,2}$,
combined with the estimates (\ref{b3}), (\ref{b5}). Similar arguments based on the Sobolev inequality and (\ref{b4}), (\ref{b5})
yield (\ref{b7}).

\subsection{Convergence}

To begin, we denote
\begin{equation}\label{l4}
\alpha = \frac{1}{\Ov{\vr}} \frac{\partial p (\Ov{\vr}, \Ov{\vt})}{\partial \vr}, \
\beta = \frac{1}{\Ov{\vr}} \frac{\partial p (\Ov{\vr}, \Ov{\vt})}{\partial \vt},\
\delta = \Ov{\vr} \frac{\partial s (\Ov{\vr}, \Ov{\vt})}{\partial \vt}, \ a(\Ov{\vr}, \Ov{\vt}) = \frac{1}{\Ov{\vr}} \frac{\beta}{\alpha}.
\end{equation}
It follows from (\ref{b4}--\ref{b5}) and the structural assumptions on the pressure (\ref{i10}), (\ref{i13}--\ref{i15}) that
\bFormula{Z1-}
\Big[\frac{\vre-\overline\vr_\ep}\ep\Big]_{\rm res} \to 0 \ \mbox{in}\ L^\infty(0,T; L^{5/3}(\Omega)),\
\Big[{ \frac{\vte-\overline\vt}\ep}\Big]_{\rm res} \to 0 \ \mbox{in}\ L^\infty(0,T; L^{4}(\Omega)).
\eF

Next, writing
\[
\frac{1}{\ep} \Grad p(\vre, \vte) - \vre \Grad F = \frac{1}{\ep} \Grad p(\vre, \vte) - \Ov{\vr}_\ep \Grad F + \ep \frac{\Ov{\vr}_\ep - \vre}{\ep} \Grad F
= \frac{1}{\ep} \Grad \left( p(\vre, \vte) - p(\Ov{\vr}_\ep, \Ov{\vt}) \right) + \ep \frac{\Ov{\vr}_\ep - \vre}{\ep} \Grad F,
\]
we deduce from the momentum balance (\ref{v2}) that
\begin{equation}\label{Z1}
\alpha\Big[\frac{\vre-\overline\vr_\ep}\ep\Big]_{\rm ess} +
\beta\Big[\frac{\vte-\overline\vt}\ep\Big]_{\rm ess} \to 0 \ \mbox{weakly-(*) in}\ L^\infty (0,T; L^2(\Omega)).
\end{equation}

Finally, we use (\ref{b3}--\ref{b5}) to show that
\bFormula{conv1}
\vre \vue \to \Ov{ \vr \vu } \ \mbox{weakly-(*) in}\ L^\infty(0,T; L^2 + L^{5/4}(\Omega;R^3)),
\eF
where, passing to the limit in the continuity equation (\ref{v1}), we may infer that
\bFormula{conv2}
\Div (\Ov{ \vr \vu }) = 0.
\eF

%
%

\section{Acoustic and thermal energy transport equations}
\label{r}

Similarly to \cite{FeiNov12}, our aim is to use the relative entropy inequality (\ref{r1}) to deduce the convergence to the target system.
To this end, we take
\[
\vr = \vre,\ \vt = \vte, \ \vu = \vue
\]
and choose the test functions $\{ r, \Theta, \vc{U} \}$ in the
following way:
\bFormula{Ap}
r = r_\ep = \Ov{\vr}_\ep + \ep R_\ep ,\
\Theta = \Theta_\ep = \Ov{\vt} + \ep T_\ep , \ \vc{U} = \vc{U}_\ep
= \vc{v} + \Grad \Phi_\ep; \eF
where $\vc{v}$ is the velocity component of the solution to
the incompressible Euler-Boussinesq system (\ref{l1})-(\ref{l3}), with the
initial condition (\ref{m6}), and the functions $R_\ep$, $T_\ep$, and $\Phi_\ep$
satisfy the \emph{acoustic equation}:
\bFormula{r3} \ep \partial_t
(\alpha R_\ep + \beta T_\ep ) +  \omega \Delta \Phi_\ep = 0, \eF
\bFormula{r4} \ep \partial_t \Grad \Phi_\ep + \Grad (\alpha R_\ep
+ \beta T_\ep ) = 0, \ \Grad \Phi_\ep \cdot \vc{n}|_{\partial \Omega} = 0,
\eF
with the initial values determined by
\bFormula{r4b} R_\ep (0, \cdot) = R_{0},\ T_\ep (0, \cdot) =
T_{0},\ \Phi_{\ep}(0, \cdot) = \Phi_{0}, \eF
and
the constants $\alpha$, $\beta$ defined in (\ref{l4}),
\[
\omega = \Ov{\vr} \left( \alpha + \frac{\beta^2}{\delta} \right).
\]

The first equation in (\ref{r3}) is nothing other than a linearization of the continuity equation, while the second
equation is a linearization of the momentum equation projected onto the space of gradients.

In order to determine $R_\ep$ and $T_\ep$ in a unique way, we require
$\delta R_\ep-\beta T_\ep$, with $\delta$ defined in (\ref{l4}),
to satisfy the \emph{transport equation}
\bFormula{r4a}
\partial_t (\delta T_\ep - \beta R_\ep) + \vc{U}_\ep \cdot \Grad \left( \delta T_\ep - \beta R_\ep - \frac{\beta}{\alpha} F \right)  = 0,
\eF
where the initial data are determined by (\ref{r4b}). Equation (\ref{r4a}) is obviously related to the limit equation (\ref{l3}).
Observe that the system of linear equations (\ref{r3}--\ref{r4a}) is well-posed.

\subsection{Initial data}

In view of the future application of the relative entropy inequality (\ref{r1}), the initial data for the test functions must be taken is such a way that
\bFormula{data}
\vc{v}(0, \cdot) = \vc{v}_0 = \vc{H}[\vu_0],\
\Phi_\ep (0, \cdot) = \Phi_{0, \eta}, \ \Grad \Phi_{0, \eta} \to \vc{H}^\perp [\vu_0] \ \mbox{in}\ L^2(\Omega;R^3)\  \mbox{as}\ \eta \to 0,
\eF
\bFormula{data1}
R_\ep (0, \cdot) = R_{0, \eta}, \ \| R_{0, \eta} \|_{L^\infty(\Omega)} < c(\eta),\ R_{0, \eta} \to \vr^{(1)}_0 \ \mbox{in}\ L^2(\Omega)
\ \mbox{as}\ \eta \to 0,
\eF
and
\bFormula{data2}
T_\ep (0, \cdot) = T_{0, \eta}, \ \| T_{0, \eta} \|_{L^\infty(\Omega)} < c(\eta),\ T_{0, \eta} \to \vt^{(1)}_0 \ \mbox{in}\ L^2(\Omega)
\ \mbox{as}\ \eta \to 0.
\eF
Note that (\ref{data} - \ref{data2}) imply that
\bFormula{data3}
\mathcal{E}_\ep \left( \vr_{0,\ep}, \vt_{0, \ep}, \vu_{0,\ep} \Big| r_\ep(0,\cdot), \Theta_\ep (0, \cdot), \vc{U}(0, \cdot) \right)
\to \chi(\eta) \ \mbox{as} \ \ep \to 0,
\eF
where
\[
\chi(\eta) \to 0 \ \mbox{as}\ \eta \to 0.
\]

Our next goal is to choose suitable approximations for the initial data.
Following \cite{FeNoSun2}, we consider the Neumann Laplacean $\Delta_N$,
\[
\mathcal{D}(\Delta_N) = \left\{ v \in L^2(\Omega) \ \Big| \ \Grad v \in L^2(\Omega;R^3), \ \intO{ \Grad v \cdot \Grad \varphi } =
\intO{ g \varphi } \right.
\]
\[
\left.
 \mbox{for any}\ \varphi \in \DC(\Ov{\Omega}) \ \mbox{and a certain}\ g \in L^2(\Omega) \right\},
\]
together with a family of regularizing operators
\bFormula{c1--}
[v]_{\eta} = G_\eta(\sqrt{- \Delta_N}) [ \psi_{1/ \eta} v ],
\eF
with the cut-off functions
\begin{equation}\label{cutoff}
\psi_\eta(x)=\psi(x/\eta);\; \psi\in C^\infty_c(R),\;0\le\psi\le 1,\; \psi(x)=
\left\{\begin{array}{c}
1\;\mbox{si $|x|\le 1$, }\\
0\;\mbox{si $|x|\ge 2$}
\end{array}\right\},
\end{equation}
\[
G_\eta \in \DC(R),\
0 \leq G_\eta \leq 1, \ G_\eta (-z) = G_\eta (z),
\]
\[
G_\eta (z) = 1 \ \mbox{for}\ z \in \left( - \frac{1}{\eta}, - \eta \right) \cup \left( {\eta},
\frac{1}{\eta}
\right), \ G_\eta(z) = 0 \ \mbox{for}\ z \in \left(-\infty, - \frac{2}{\eta} \right) \cup
\left(- \frac{\eta}{2}, \frac{\eta}{2} \right) \cup \left( \frac{2}{ \eta}, \infty \right),
\]
where the linear operator $G_\eta(\sqrt{- \Delta_N})$ is defined by means of the standard spectral theory associated to $\Delta_N$.

Accordingly, we consider regularized initial data in the form
\bFormula{d1}
R_{0, \eta} = [ \vr^{(1)}_{0} ]_{\eta} ,\ T_{0, \eta} = [ \vt^{(1)}_{0}]_\eta,
\eF
and
\bFormula{d2}
\Phi_{0, \eta} = \Big[ \Delta^{-1}_N \Div [\vu_{0}] \Big]_\eta , \ \mbox{with}\
\Grad \Delta^{-1}_N \Div [\vu_{0}] \equiv \vc{H}^\perp [ \vu_{0} ].
\eF

To avoid excessive notation, we omit writing the parameter $\eta$ in the course of the limit passage $\ep \to 0$.

\subsection{Dispersive estimates for the wave equation}

The acoustic equation (\ref{r3} - \ref{r4b}) has been studied in detail in \cite{FeNoSun2}. In particular,
we report the following estimates
(\cite[estimates (6.6), (6.8)]{FeNoSun2}:
\bFormula{a1}
\sup_{t \in [0,T]} \left( \left\| \Grad \Phi_{\ep, \eta} \right\|_{W^{k,2} \cap W^{k,\infty}(\Omega:R^3)} +
\left\| (\alpha R_{\ep, \eta} + \beta T_{\ep, \eta})  (t,\cdot) \right\|_{W^{k,2} \cap W^{k,\infty}(\Omega:R^3)} \right)
\eF
\[
\leq c(k, \eta) \left(
\left\| \Grad \Phi_{0, \eta } \right\|_{L^2(\Omega;R^3)} + \left\| \alpha R_{0,\eta} + \beta T_{0, \eta} \right\|_{L^2(\Omega)} \right),\
\]
for any $k=0,1,\dots$, $\eta > 0$; and the dispersive estimates
\bFormula{a2}
\int_0^T \left( \left\| \Grad \Phi_{\ep, \eta} \right\|_{W^{k,\infty}(\Omega:R^3)} +
\left\| (\alpha R_{\ep, \eta} + \beta T_{\ep, \eta})  (t,\cdot) \right\|_{W^{k,\infty}(\Omega:R^3)} \right) \ \dt
\eF
\[
\leq \omega(\ep, \eta, k) \left(
\left\| \Grad \Phi_{0,\eta} \right\|_{L^2(\Omega;R^3)} + \left\| \alpha R_{0,\eta} + \beta T_{0, \eta}  \right\|_{L^2(\Omega)} \right)
\]
where
\[
\omega(\ep, \eta,k) \to 0 \ \mbox{as}\ \ep \to 0 \ \mbox{for any fixed}\ \eta > 0, \ k \geq 0.
\]
The relation (\ref{a2}) represents \emph{dispersive} estimates for the wave equation (\ref{r3}), (\ref{r4}). Note that both (\ref{a1}) and (\ref{a2})
apply to the regularized initial data, meaning for a fixed $\eta > 0$; they in fact blow up when $\eta \to 0$.

Moreover, as shown in \cite[Section 5.3]{EF100},
\bFormula{c1b--}
|x|^s |\partial^k_x [h]_\eta (x) | \leq c(s,k, \eta) \| h \|_{L^2(\Omega)} \ \mbox{for all}\
x \in \Omega, \ s \geq 0, k \geq 0,
\eF
therefore the functions $\Phi_{\ep, \eta}$, $(\alpha R_{\ep, \eta} + \beta T_{\ep, \eta}) $ decay fast for $|x| \to \infty$ as long as $\eta > 0$ is fixed.

\bRemark{r11}
As a matter of fact, the results of \cite{FeNoSun2} are stated for the domain $\Omega$ - a perturbed half-space. However, as pointed out in
\cite{FeNoSun2}, the same holds for a larger class
of domains on which $\Delta_N$, among which the exterior domains in $R^3$. Alternatively, we may also use the dispersive estimates established
by Isozaki \cite{Isoz}.
\eR
\subsection{$L^p$ estimates for the transport equation}

For fixed $\eta > 0$, the initial data for the transport equation (\ref{r4a}) enjoy the decay properties (\ref{c1b--}).
Consequently, in view of (\ref{a1}), (\ref{a2}), the solutions of the transport equation (\ref{r4a}) admit the estimates
\bFormula{a3} \sup_{t \in [0,T]} \left\| \delta
T_{\ep, \eta} - \beta R_{\ep, \eta} \right\|_{W^{k,q}(\Omega)} \leq c(\eta, k,F) \left(1 +
\left\| \delta T_{0,\eta} - \beta R_{0, \eta}
\right\|_{L^2(\Omega)} \right),\ k=0,1,\ 1 \leq q \leq \infty, \eF
and the family
\bFormula{a4}
\left\{ \delta
T_{\ep, \eta} - \beta R_{\ep, \eta} \right\}_{\ep > 0} \ \mbox{is precompact in} \ C([0,T]; W^{k,q}(\Omega)), \ k=0,1,\ 1 \leq q \leq \infty.
\eF

Consequently, combining (\ref{a2}), (\ref{a3}), (\ref{a4}) we can let $\ep \to 0$ to obtain
\begin{equation}\label{a5}
T_{\ep, \eta}  \to\ T_\eta \;\mbox{ strongly in $L^\infty_{\rm loc}((0,T]; W^{k,p}(\Omega))$, \ $p>2$, and weakly$-(*)$ in $L^\infty(0,T; W^{k,2}(\Omega))$,
\ $k=0,1$,}
\ \mbox{as}\ \ep \to 0,
\end{equation}
\begin{equation}\label{a5aa}
R_{\ep, \eta}  \to\ R_\eta \;\mbox{ strongly in $L^\infty_{\rm loc}((0,T]; W^{k,p}(\Omega))$, \ $p>2$, and weakly$-(*)$ in
$L^\infty(0,T; W^{k,2}(\Omega))$,\ $k=0,1$,}
\ \mbox{as}\ \ep \to 0,
\end{equation}
where $T_\eta$ satisfies
\bFormula{a5A}
c_p(\Ov{\vr}, \Ov{\vt}) \left( \partial_t T_\eta + \vc{v} \cdot \Grad T_\eta \right) - \Ov{\vt} a(\Ov{\vr}, \Ov{\vt}) \vc{v} \cdot \Grad F = 0,
\eF
with the initial data
\bFormula{a5B}
T_\eta (0, \cdot) = \frac{\Ov{\vt}}{c_p(\Ov{\vr}, \Ov{\vt})} \left( \frac{\partial s(\Ov{\vr}, \Ov{\vt})}{
\partial \vr} [\vr^{(1)}_0]_{\eta} + \frac{\partial s(\Ov{\vr}, \Ov{\vt})}{
\partial \vt} [\vt^{(1)}_0]_\eta \right).
\eF

\section{Convergence}
\label{c}

In this section, we use the test functions
(\ref{Ap}) in the relative entropy inequality (\ref{r1}). Fixing $\eta > 0$ we perform the limit for $\ep \to 0$.
This will be carried over in several steps in the spirit of \cite{FeiNov12}. We omit the subscript $\eta$ whenever no confusion arises.

\subsection{{ Viscous and heat conducting terms}}

We show by direct calculation, splitting the terms in their essential and residual parts
and using  assumptions (\ref{i8}--\ref{i9}), uniform bounds (\ref{b5}--\ref{b7}), regularity (\ref{euler}), and  estimates (\ref{a1}--\ref{a4})
that the dissipative terms related to the viscosity and to
the heat conductivity on the right-hand side of (\ref{r1}) become
negligible as $\ep \to 0$. More precisely:
\[
\ep^a \tn{S}(\vte, \Grad \vue) : \Grad \vc{U}_\ep \to 0 \ \mbox{in} \ L^2((0,T) \times \Omega)+
 L^2(0,T; L^{4/3}(\Omega;R^3)) \ \mbox{as}\ \ep
\to 0,
\]
and
\[
\ep^{b - 2} \frac{ \vc{q}(\vte, \Grad \vte) \cdot \Grad \Theta_\ep }{\vte}  \to 0
\ \mbox{in}\ L^2((0,T) \times \Omega) + L^1((0,T) \times \Omega) \ \mbox{as}\ \ep
\to 0.
\]

Consequently, combining the previous observation with (\ref{data3}), we can write the relative entropy inequality (\ref{r1})
as

\bFormula{r9}
\mathcal{E}_\ep \left( \vre, \vte, \vue \Big| r_\ep , \Theta_\ep , \vc{U}_\ep \right) (\tau)
\eF
\[
\leq \chi (\ep, \eta)  + \int_0^\tau \intO{  \vre \Big( \partial_t \vc{U}_\ep + \vue
\cdot\Grad \vc{U}_\ep \Big)\cdot (\vc{U}_\ep - \vue)    } \ \dt
\]
\[
- \frac{1}{\ep} \int_0^\tau \intO{ \left( \vre \Big( s(\vre,\vte) - s(r_\ep, \Theta_\ep) \Big) \partial_t T_\ep + \vre \Big( s(\vre,\vte) - s(r_\ep, \Theta_\ep) \Big) \vue \cdot \Grad T_\ep \right) } \ \dt
\]
$$
+\frac
1{\ep^2}\int_0^\tau\intO{\Big[\Big(p(r_\ep,\Theta_\ep)-p(\vre,\vte)\Big){\rm
div}\vc U_\ep +\frac\vre {r_\ep}(\vc U_\ep-\vue)\cdot\Grad
p(r_\ep,\Theta_\ep)\Big]}{\rm d}t
$$
\[
+ \frac{1}{\ep^2} \int_0^\tau \intO{ \frac{r_\ep-\vre}{r_\ep}\Big(
\partial_t p(r_\ep , \Theta_\ep ) +\vc U_\ep \cdot
\Grad p(r_\ep , \Theta_\ep ) \Big) } \ \dt -\frac 1\ep\int_0^\tau\int_{R^3}\vre\Grad F\cdot(\vc U_\ep-\vc u_\ep){\rm d}x{\rm d}t,
\]
where $\chi$ {  denotes a generic function} satisfying
\begin{equation}\label{chi}
{ \lim_{\eta \to 0} \Big( \lim_{\ep \to 0} \chi(\ep, \eta) \Big) }= 0.
\end{equation}

\subsection{Velocity dependent terms}

Our next goal is to handle the expression
\[
\int_0^\tau \intO{ \Big[  \vre (\vc{U}_\ep - \vue ) \cdot
\partial_t \vc{U}_\ep + \vre (\vc{U}_\ep - \vue) \otimes \vue : \Grad \vc{U}_\ep \Big] } \ \dt =
\]
\[
 \int_0^\tau \intO{ \vre (\vc{U}_\ep - \vue ) \otimes
(\vue - \vc{U}_\ep ) : \Grad \vc{U}_\ep } \ \dt
\]
\[ + \int_0^\tau \intO{ \vre (\vc{U}_\ep - \vue) \cdot \Big( \partial_t \vc{v} + \vc{v}
\cdot \Grad \vc{v} \Big) } \ \dt + \int_0^\tau \intO{ \vre (\vc{U}_\ep - \vue)
\cdot \partial_t \Grad \Phi_\ep } \ \dt
\]
\[
+ \int_0^\tau \intO{ \vre (\vc{U}_\ep - \vue ) \otimes
\Grad \Phi_\ep : \Grad \vc{v}} + \int_0^\tau \intO{ \vre (\vc{U}_\ep - \vue ) \otimes
\vc{v} : \Grad^2 \Phi_\ep } \ \dt
\]
\[
+\frac{1}{2} \int_0^\tau \intO{ \vre (\vc{U}_\ep - \vue) \cdot \Grad |\Grad \Phi_\ep |^2 } \ \dt.
\]

Thanks to (\ref{euler}), (\ref{a1}), (\ref{a2}) and the energy bounds established in (\ref{b3} - \ref{b7}), the first integral on the right hand side
can be dominated by the expression
$$
\chi (\ep, \eta) + c \int_0^\tau{\cal E}\Big(\vre,\vte, \vue \Big|r_\ep,\Theta_\ep,\vc U_\ep\Big){\rm d}t,
$$
with $c$ independent of $\ep$, $\eta$.

The second term reads
\[
\int_0^\tau \int_{\Omega}{ \vre \vue \cdot \Grad \Pi } \ \dt
- \int_0^\tau \int_{\Omega}{
\vre (\vc{v} + \Grad \Phi_\ep ) \cdot \Grad \Pi } \ \dt
\]
$$
+\frac 1{\overline\vr} \frac {\beta}{\alpha }\int_0^\tau \int_{\Omega}{ \theta\vre \vue \cdot \Grad F } \ \dt
- \frac 1{\overline\vr}\frac {\beta}{\alpha}\int_0^\tau \int_{\Omega}
\theta\vre (\vc{v} + \Grad \Phi_\ep ) \cdot \Grad F \ \dt
$$
$$
=\frac 1{\overline\vr}\frac {\beta}{\alpha}\int_0^\tau \intO{ \theta\vre \vue \cdot \Grad F } \ \dt
- \frac 1{\overline\vr}\frac {\beta}{\alpha}\int_0^\tau \int_{\Omega}
\theta\vre \vc{v}  \cdot \Grad F \ \dt
+\chi(\ep,\eta)
$$
where we have used  the equations (\ref{l1}--\ref{l2}), formulas (\ref{conv1}--\ref{conv2}), the dispersive estimates (\ref{a2}), and relation (\ref{euler}).

Next, using the equation (\ref{r4}), we may write the third integral in the form
\[
- \int_0^\tau \intO{ \vre \vue \cdot \partial_t \Grad \Phi_\ep } \ \dt
-\int_0^\tau \intO{ \frac{\vre - \Ov{\vr}}{\ep}  \vc{v} \cdot \Grad \left( \alpha R_\ep + \beta T_\ep \right) } \ \dt
\]
\[
-
\int_0^\tau \intO{ \frac{\vre - \Ov{\vr}}{\ep} \Grad \Phi_\ep \cdot \Grad (\alpha R_\ep + \beta T_\ep) } \ \dt+
\frac{1}{2} \int_0^\tau \intO{ \Ov{\vr} \partial_t | \Grad \Phi_\ep |^2 } \ \dt
\]
\[
=- \int_0^\tau \intO{ \vre \vue \cdot \partial_t \Grad \Phi_\ep } \ \dt +
\frac{1}{2} \int_0^\tau \intO{ \Ov{\vr} \partial_t | \Grad \Phi_\ep |^2 } \ \dt  +{ \chi(\ep,\eta)}.
\]
where we have used wave equation (\ref{r3}--\ref{r4}), estimates (\ref{b3}--\ref{b5}), (\ref{cutoff}), regularity of $\vc v$ stated (\ref{euler}),
the relation (\ref{p3}),
and dispersive estimates (\ref{a2}).

Finally, in view of the uniform bounds (\ref{euler}), (\ref{b3} - \ref{b5}), and the {dispersive estimates} stated in (\ref{a2}),
the last three integrals tend to zero for $\ep \to 0$, uniformly with respect to $\tau$.

Resuming, we obtain
\[
\int_0^\tau \intO{ \Big[  \vre (\vc{U}_\ep - \vue ) \cdot
\partial_t \vc{U}_\ep + \vre (\vc{U}_\ep - \vue) \otimes \vue : \Grad \vc{U}_\ep \Big] } \ \dt
\]
$$
\le \chi (\ep, \eta) + c\int_0^\tau{\cal E}\Big(\vre,\vte,\vue\Big|r_\ep,\Theta_\ep,\vc U_\ep\Big){\rm d} t +
\frac 1{\overline\vr}\frac {\beta}{\alpha}\int_0^\tau \intO{ \theta \vre\vue \cdot \Grad F } \ \dt
$$
$$
- \int_0^\tau \intO{ \vre \vue \cdot \partial_t \Grad \Phi_\ep } \ \dt +
\frac{1}{2} \int_0^\tau \intO{ \Ov{\vr} \partial_t | \Grad \Phi_\ep |^2 } \ \dt
$$
$$
- \frac 1{\overline\vr}\frac {\beta}{\alpha}\int_0^\tau \int_{\Omega}
\theta\vre \vc{v}  \cdot \Grad F \ \dt ;
$$
whence relation (\ref{r9}) becomes

\bFormula{r10}
\mathcal{E}_\ep \left( \vre, \vte, \vue \Big| r_\ep , \Theta_\ep , \vc{U}_\ep \right) (\tau) \leq \chi(\ep, \eta)  +
c \int_0^\tau  \mathcal{E}_\ep \left( \vre, \vte, \vue \Big| r_\ep , \Theta_\ep , \vc{U}_\ep \right) \ \dt
\eF
\[
+ \left[ \intO{ \Ov{\vr} \frac{1}{2}|\Grad \Phi_\ep |^2 } \right]_{t = 0}^{t = \tau} -
\int_0^\tau \intO{ \vre \vue \cdot \partial_t \Grad \Phi_\ep } \ \dt
\]
\[
- \frac{1}{\ep} \int_0^\tau \intO{ \left[ \vre \Big( s(\vre,\vte) - s(r_\ep, \Theta_\ep) \Big) \partial_t T_\ep + \vre \Big( s(\vre,\vte) - s(r_\ep, \Theta_\ep) \Big) \vue \cdot \Grad T_\ep  \right] } \ \dt
\]
\[
+ \frac{1}{\ep^2} \int_0^\tau \intO{ \left( ( r_\ep - \vre ) \frac{1}{r_\ep } \partial_t p(r_\ep , \Theta_\ep ) - \frac{\vre}{r_\ep} \vue \cdot
\Grad p(r_\ep , \Theta_\ep ) \right) } \ \dt - \frac{1}{\ep^2} \int_0^\tau \intO{ \Big( p(\vre, \vte) - p(\Ov{\vr}_\ep, \Ov{\vt}) \Big)
\Delta \Phi_\ep } \ \dt
\]
\[
-\frac 1\ep\int_0^\tau\int_{\Omega}\vre\Grad F\cdot(\vc{v} -\vc u_\ep) \ \dxdt
+
\frac 1{\overline\vr}\frac {\beta}{\alpha}\int_0^\tau \int_{\Omega}{ \theta\vre \vue \cdot \Grad F }{\rm d} x \ \dt
- \frac 1{\overline\vr}\frac {\beta}{\alpha}\int_0^\tau \int_{\Omega}
\theta\vre \vc{v}  \cdot \Grad F \ \dxdt.
$$
In the above, we have used the identity
\[
\intO{ \left[ \Big( p(r_\ep, \Theta_\ep) - p(\vre, \vte) \Big) \Div \vc{U}_\ep + \left( 1 - \frac{\vre}{r_\ep} \right) \vc{U}_\ep
\cdot \Grad p(r_\ep, \Theta_\ep) + \frac{\vre}{r_\ep} (\vc{U}_\ep - \vue) \cdot \Grad p(r_\ep, \Theta_\ep)  \right] }
\]
\[
= - \intO{  p(\vre, \vte)  \Delta \Phi_\ep}  - \intO{ \frac{\vre}{r_\ep} \vue \cdot \Grad p(r_\ep, \Theta_\ep) },
\]
together with
\[
-\frac 1\ep\int_0^\tau\int_{\Omega}\vre\Grad F\cdot(\vc U_\ep-\vc u_\ep) \ \dxdt =
-\frac 1\ep\int_0^\tau\int_{\Omega}\vre\Grad F\cdot(\vc{v}-\vc u_\ep) \ \dxdt - \frac 1\ep\int_0^\tau\int_{\Omega}\vre\Grad F\cdot\Grad \Phi_\ep  \ \dxdt
\]
\[
= { \chi(\ep,  \eta)} -\frac 1\ep\int_0^\tau\int_{\Omega}\vre\Grad F\cdot(\vc{v}-\vc u_\ep) \ \dxdt + \frac{1}{\ep^2} p(\Ov{\vr}_\ep, \Ov{\vt}) \Delta \Phi_\ep
\ \dxdt.
\]

Recall that $\Grad \Phi_\ep(t, \cdot)$ decays fast as $|x| \to \infty$ and $\Div \vc{v} = 0$, which justifies the by-parts integration.

\subsection{Pressure dependent terms}

We write
\[
\frac{1}{\ep^2}  \frac{\vre}{r_\ep} \vue \cdot \Grad p(r_\ep,
\Theta_\ep) = \frac{1}{\ep^2}  \frac{\vre}{r_\ep} \vue \cdot \Grad \Big( p(r_\ep, \Theta_\ep) - p(\Ov{\vr}_\ep, \Ov{\vt}) \Big) +
\frac{1}{\ep^2}  \frac{\vre}{r_\ep} \vue \cdot \Grad p(\Ov{\vr}_\ep, \Ov{\vt})
\]
\[
= \frac{1}{\ep^2}  \frac{\vre}{r_\ep} \vue \cdot \Grad \left( p(r_\ep, \Theta_\ep) - \frac{\partial p(\Ov{\vr}_\ep, \Ov{\vt}) }
{\partial \vr} \ep R_\ep - \frac{\partial p(\Ov{\vr}_\ep, \Ov{\vt}) }
{\partial \vt} \ep T_\ep  - p(\Ov{\vr}_\ep, \Ov{\vt}) \right)
\]
\[
+\frac{1}{\ep}  \frac{\vre}{r_\ep} \vue \cdot \Grad \left( \frac{\partial p(\Ov{\vr}_\ep, \Ov{\vt}) }
{\partial \vr} R_\ep + \frac{\partial p(\Ov{\vr}_\ep, \Ov{\vt}) }
{\partial \vt} T_\ep   \right) + \frac{1}{\ep} \frac{\Ov{\vr}_\ep} {r_\ep} \vre \vue \cdot \Grad F.
\]

Next, we use the decay properties of the equilibrium density profile $\Ov{\vr}_\ep$ stated in (\ref{p3}), together with (\ref{a5}), (\ref{a5aa})
to observe that
\[
\frac{1}{\ep^2 r_\ep} \Grad
\left( p(r_\ep, \Theta_\ep) - \frac{\partial p(\Ov{\vr}_\ep, \Ov{\vt}) }
{\partial \vr} \ep R_\ep - \frac{\partial p(\Ov{\vr}_\ep, \Ov{\vt}) }
{\partial \vt} \ep T_\ep  - p(\Ov{\vr}_\ep, \Ov{\vt}) \right) \to \Grad H \ \mbox{in}\ L^p (0,T; (L^2 \cap L^q)(\Omega;R^3)) ,\ p \geq 1,\ q > 2,
\]
where the right-hand side is a gradient of a certain function $H$. Consequently, using (\ref{conv1}), (\ref{conv2}) we may infer that
\[
\int_0^\tau \intO{ \frac{1}{\ep^2}  \frac{\vre}{r_\ep} \vue \cdot \Grad \left( p(r_\ep, \Theta_\ep) - \frac{\partial p(\Ov{\vr}_\ep, \Ov{\vt}) }
{\partial \vr} \ep R_\ep - \frac{\partial p(\Ov{\vr}_\ep, \Ov{\vt}) }
{\partial \vt} \ep T_\ep  - p(\Ov{\vr}_\ep, \Ov{\vt}) \right) } \ \dt = { \chi (\ep, \eta)}.
\]

Moreover, by the same token, we obtain
\[
\int_0^\tau \intO{
\frac{1}{\ep}  \frac{\vre}{r_\ep} \vue \cdot \Grad \left( \frac{\partial p(\Ov{\vr}_\ep, \Ov{\vt}) }
{\partial \vr} R_\ep + \frac{\partial p(\Ov{\vr}_\ep, \Ov{\vt}) }
{\partial \vt} T_\ep   \right) } \ \dt = \eta(\ep, \delta) + \int_0^\tau \intO{
\frac{1}{\ep}  {\vre} \vue \cdot \Grad \left( \alpha
 R_\ep + \beta
 T_\ep   \right) } \ \dt.
\]

Making use of the identity
\[
\int_0^\tau \intO{
\frac{1}{\ep}  {\vre} \vue \cdot \Grad \left( \alpha
 R_\ep + \beta
 T_\ep   \right) } \ \dt = - \int_0^\tau\int_{R^3}\vre \vue \cdot \partial_t \Grad \Phi_\ep{\rm d} x{\rm d} t
\]
we may rewrite (\ref{r10}) in the form
\bFormula{r11-}
\mathcal{E}_\ep \left( \vre, \vte, \vue \Big| r_\ep , \Theta_\ep , \vc{U}_\ep \right) (\tau) \leq \chi(\ep, \eta)  +
c \int_0^\tau  \mathcal{E}_\ep \left( \vre, \vte, \vue \Big| r_\ep , \Theta_\ep , \vc{U}_\ep \right) \ \dt
+ \left[ \intO{ \Ov{\vr} \frac{1}{2}|\Grad \Phi_\ep |^2 } \right]_{t = 0}^{t = \tau}
\eF
\[
- \frac{1}{\ep} \int_0^\tau \intO{ \left[ \vre \Big( s(\vre,\vte) - s(r_\ep, \Theta_\ep) \Big) \partial_t T_\ep + \vre \Big( s(\vre,\vte) - s(r_\ep, \Theta_\ep) \Big) \vue \cdot \Grad T_\ep  \right] } \ \dt
\]
\[
+ \frac{1}{\ep^2} \int_0^\tau \intO{  ( r_\ep - \vre ) \frac{1}{r_\ep } \partial_t p(r_\ep , \Theta_\ep )  } \ \dt - \frac{1}{\ep^2} \int_0^\tau \intO{ \Big( p(\vre, \vte) - p(\Ov{\vr}_\ep, \Ov{\vt}) \Big)
\Delta \Phi_\ep } \ \dt
\]
\[
+ \int_0^\tau\int_{\Omega}\frac{R_\ep}{r_\ep} \vre \vue \cdot \Grad F  \ \dxdt - \int_0^\tau \intO{ \frac{\vre - \Ov{\vr}_\ep }{\ep} \vc{v} \cdot \Grad F
} \ \dt
\]
\[
+
\frac 1{\overline\vr}\frac {\beta}{\alpha}\int_0^\tau \int_{\Omega}{ \theta\vre \vue \cdot \Grad F }{\rm d} x \ \dt
- \frac 1{\overline\vr}\frac {\beta}{\alpha}\int_0^\tau \int_{\Omega}
\theta\vre \vc{v}  \cdot \Grad F \ \dxdt.
\]

Finally, we use the fact that
\bFormula{brum}
\alpha R_{\eta} + \beta T_{\eta} = 0,
\eF
and that $T_\eta$ and $\theta$ satisfy the same equation { (see (\ref{a5A}) and (\ref{l3}))}  with the initial data given by
(\ref{a5B}), (\ref{m6}), respectively, to deduce that

\bFormula{r11}
\mathcal{E}_\ep \left( \vre, \vte, \vue \Big| r_\ep , \Theta_\ep , \vc{U}_\ep \right) (\tau) \leq \chi(\ep, \eta)  +
c \int_0^\tau  \mathcal{E}_\ep \left( \vre, \vte, \vue \Big| r_\ep , \Theta_\ep , \vc{U}_\ep \right) \ \dt
+ \left[ \intO{ \Ov{\vr} \frac{1}{2}|\Grad \Phi_\ep |^2 } \right]_{t = 0}^{t = \tau}
\eF
\[
- \frac{1}{\ep} \int_0^\tau \intO{ \left[ \vre \Big( s(\vre,\vte) - s(r_\ep, \Theta_\ep) \Big) \partial_t T_\ep + \vre \Big( s(\vre,\vte) - s(r_\ep, \Theta_\ep) \Big) \vue \cdot \Grad T_\ep  \right] } \ \dt
\]
\[
+ \frac{1}{\ep^2} \int_0^\tau \intO{  ( r_\ep - \vre ) \frac{1}{r_\ep } \partial_t p(r_\ep , \Theta_\ep )  } \ \dt - \frac{1}{\ep^2} \int_0^\tau \intO{ \Big( p(\vre, \vte) - p(\Ov{\vr}_\ep, \Ov{\vt}) \Big)
\Delta \Phi_\ep } \ \dt
\]
\[
 - \int_0^\tau \intO{ \frac{\vre - \Ov{\vr}_\ep }{\ep} \vc{v} \cdot \Grad F
} \ \dt
- \frac 1{\overline\vr}\frac {\beta}{\alpha}\int_0^\tau \int_{\Omega}
\theta\vre \vc{v}  \cdot \Grad F \ \dxdt.
\]

\subsection{Replacing velocity in the entropy convective term}\label{7.4}

Our intention in this section is to ``replace'' $\vue$ by $\vc{U}_\ep$ in the remaining (last) convective term in (\ref{r11}). To this end,
we write
\[
\int_0^\tau \intO{ \vre \frac{ s(\vre, \vte) - s(r_\ep, \Theta_\ep ) }{ \ep } \vue \cdot \Grad T_\ep } \ \dt
\]
\[
=
\int_0^\tau \intO{ \vre \frac{ s(\vre, \vte) - s(r_\ep , \Theta_\ep ) }{ \ep } \vc{U}_\ep  \cdot \Grad T_\ep } \ \dt +
\int_0^\tau \intO{ \vre \frac{ s(\vre, \vte) - s(r_\ep, \Theta_\ep ) }{ \ep } (\vue - \vc{U}_\ep)  \cdot \Grad T_\ep } \ \dt,
\]
where
\[
\left| \int_0^\tau \intO{ \vre \left[ \frac{ s(\vre, \vte) - s(r_\ep , \Theta_\ep ) }{ \ep } \right]_{\rm ess} (\vue - \vc{U}_\ep)  \cdot \Grad T_\ep } \ \dt \right|
\]
\[
\leq A(\eta) \int_0^\tau \intO{ \left( \vre | \vue - \vc{U}_\ep |^2
+ \left| \left[ \frac{\vre - r_\ep}{\ep} \right]_{\rm ess} \right|^2 + \left| \left[ \frac{\vte - \Theta_\ep}{\ep} \right]_{\rm ess} \right|^2 \right)
}\ \dt
\]
$$
\le c \int_0^\tau{\cal E}\Big(\vre,\vte,\vue\Big|r_\ep,\Theta_\ep,\vc U_\ep\Big){\rm d}t
$$
and
\[
+ \int_0^\tau \intO{ \vre \left[ \frac{ s(\vre, \vte) - s(r_\ep , \Theta_\ep ) }{ \ep } \right]_{\rm res} (\vue - \vc{U}_\ep)  \cdot \Grad T_\ep } \
\dt=\chi(\ep ,\eta)\;\mbox{provided $0<a<10/3$.}
\]
When estimating the residual component, we have first deduced from (\ref{i12} - \ref{i16}) the inequality
\begin{equation}\label{*entropy}
\vr | s(\vr, \vt) | \leq c\left( \vt^3 + \vr |\log(\vr)| + \vr [\log(\vt)]^+ \right)
\end{equation}
and then employed the estimates (\ref{b5}--\ref{b6}) for $\vre$, $\vte$, together with the estimates (\ref{a2}--\ref{a4}) for $R_\ep$, $T_\ep$,
$\Grad\Phi_\ep$, and (\ref{euler}) for $\vc v$.

Consequently, we can can rewrite inequality (\ref{r11}) in the form
\bFormula{r11+}
\mathcal{E}_\ep \left( \vre, \vte, \vue \Big| r_\ep , \Theta_\ep , \vc{U}_\ep \right) (\tau) \leq \chi(\ep, \eta)  +
c \int_0^\tau  \mathcal{E}_\ep \left( \vre, \vte, \vue \Big| r_\ep , \Theta_\ep , \vc{U}_\ep \right) \ \dt
+ \left[ \intO{ \Ov{\vr} \frac{1}{2}|\Grad \Phi_\ep |^2 } \right]_{t = 0}^{t = \tau}
\eF
\[
- \frac{1}{\ep} \int_0^\tau \intO{ \left[ \vre \Big( s(\vre,\vte) - s(r_\ep, \Theta_\ep) \Big) \partial_t T_\ep + \vre \Big( s(\vre,\vte) - s(r_\ep, \Theta_\ep)
\Big) \vc{U}_\ep \cdot \Grad T_\ep  \right] } \ \dt
\]
\[
+ \frac{1}{\ep^2} \int_0^\tau \intO{  ( r_\ep - \vre ) \frac{1}{r_\ep } \partial_t p(r_\ep , \Theta_\ep )  } \ \dt - \frac{1}{\ep^2} \int_0^\tau \intO{ \Big( p(\vre, \vte) - p(\Ov{\vr}_\ep, \Ov{\vt}) \Big)
\Delta \Phi_\ep } \ \dt
\]
\[
 - \int_0^\tau \intO{ \frac{\vre - \Ov{\vr}_\ep }{\ep} \vc{v} \cdot \Grad F
} \ \dt
- \frac 1{\overline\vr}\frac {\beta}{\alpha}\int_0^\tau \int_{\Omega}
\theta\vre \vc{v}  \cdot \Grad F \ \dxdt.
\]

\subsection{The entropy and the pressure}
{ \subsubsection{Handling the residual component}}
To begin, we observe
that the residual components of all integrals on the second and third line of inequality (\ref{r11+}) are negligible.
To this end, we first use the estimates (\ref{a2} - \ref{a4}), (\ref{a5}), (\ref{a5aa}), together with the equations (\ref{r3} - \ref{r4a}),
to deduce
\bFormula{pom1}
\sup_{t \in [0,T]} \ep \| \partial_t R_\ep (t, \cdot) \|_{L^\infty(R^3)} ,\
\sup_{t \in [0,T]} \ep \| \partial_t T_\ep (t, \cdot) \|_{L^\infty(R^3)} \leq A(\eta),
\eF
\bFormula{pom2}
\ep \| \partial_t R_\ep (t, \cdot) \|_{L^\infty(R^3)} \to 0 , \ \ep \| \partial_t T_\ep (t, \cdot) \|_{L^\infty(R^3)} \to 0\ \mbox{for any} \ t > 0.
\eF
Now, we employ these relations in combination with the uniform estimates (\ref{b5}); after a long but straightforward calculation,
we finally get the desired result, namely
\begin{equation}\label{*1-}
- \frac{1}{\ep} \int_0^\tau \intO{ \left[ \Big[\vre \Big( s(\vre,\vte)
- s(r_\ep, \Theta_\ep) \Big) \partial_t T_\ep + \vre \Big(
s(\vre,\vte) - s(r_\ep, \Theta_\ep) \Big) \vc U_\ep \cdot \Grad
T_\ep\Big]_{\rm res} \right] } \ \dt
\end{equation}
$$
- \frac{1}{\ep^2} \int_0^\tau \intO{ \Big[\frac{\vr_\ep - r_\ep}
{r_\ep}\partial_t p(r_\ep , \Theta_\ep ) \Big]_{\rm res} } \ \dt -\frac 1{\ep^2}
\int_0^\tau
\intO{\Big[\Big(p(\vre,\vte)-p(\overline\vr_\ep,\overline\vt)\Big)\Delta\Phi_\ep\Big]_{\rm res}}\ {\rm
d}t=\chi(\ep,\eta)
$$
\subsubsection{Handling the essential component}
In view of the preceding Section, we have to handle solely the essential part of the integrals at the first and second line of
formula (\ref{r11+}) whose integrands can be, roughly speaking, replaced by their linearization
at $\overline\vr_\ep$, $\overline\vt$. Since we already know that the corresponding residual components are negligible,
we may omit the symbol $[\cdot]_{\rm ess}$ in all integrands.

We check that
\begin{equation}\label{*1}
- \frac{1}{\ep} \int_0^\tau \intO{ \left[ \vre \Big( s(\vre,\vte)
- s(r_\ep, \Theta_\ep) \Big) \partial_t T_\ep + \vre \Big(
s(\vre,\vte) - s(r_\ep, \Theta_\ep) \Big) \vc U_\ep \cdot \Grad
T_\ep \right] } \ \dt
\end{equation}
$$
- \frac{1}{\ep^2} \int_0^\tau \intO{ \frac{\vr_\ep - r_\ep}
{r_\ep}\partial_t p(r_\ep , \Theta_\ep )  } \ \dt -\frac 1{\ep^2}
\int_0^\tau
\intO{\Big(p(\vre,\vte)-p(\overline\vr_\ep,\overline\vt)\Big)\Delta\Phi_\ep}\ {\rm
d}t
$$
$$
= -\int_0^\tau\intO{\Big(\delta\frac{\vte-\Theta_\ep}\ep
-\beta\frac{\vr_\ep-r_\ep}\ep\Big)\Big(\partial_t T_\ep + \vc
U_\ep\cdot\Grad T_\ep\Big)}\ {\rm d}t
$$
$$
-\int_0^\tau\intO{\frac{\vre- r_\ep}\ep\partial_t\Big(\alpha
R_\ep+\beta T_\ep\Big)}\ {\rm d}t +\int_0^\tau
\intO{\frac{\delta}{\beta^2+\alpha\delta}\Big(\alpha\frac{\vre
-\overline\vr_\ep}\ep+\beta\frac{\vte-\overline\vt}\ep\Big)
\partial_t\Big(\alpha R_\ep+\beta T_\ep\Big)}\ {\rm d}t
$$
\[
+ \int_0^\tau \intO{ \frac{1}{\ep} \left( \frac{\partial p(\Ov{\vr}, \Ov{\vt}) }{\partial \vr} -
\frac{\partial p(\Ov{\vr}_\ep , \Ov{\vt}) }{\partial \vr} \right) \frac{\vre - \Ov{\vr}_\ep }{\ep} \Delta \Phi_\ep } \ \dt
+ \int_0^\tau \intO{ \frac{1}{\ep} \left( \frac{\partial p(\Ov{\vr}, \Ov{\vt}) }{\partial \vt} -
\frac{\partial p(\Ov{\vr}_\ep , \Ov{\vt}) }{\partial \vt} \right) \frac{\vte - \Ov{\vt} }{\ep} \Delta \Phi_\ep } \ \dt
+ \chi(\ep, \eta),
\]
where, in accordance with the dispersive estimates (\ref{a1}), (\ref{a2}) and (\ref{p3}),
\[
\int_0^\tau \intO{ \frac{1}{\ep} \left( \frac{\partial p(\Ov{\vr}, \Ov{\vt}) }{\partial \vr} -
\frac{\partial p(\Ov{\vr}_\ep , \Ov{\vt}) }{\partial \vr} \right) \frac{\vre - \Ov{\vr}_\ep }{\ep} \Delta \Phi_\ep } \ \dt
+ \int_0^\tau \intO{ \frac{1}{\ep} \left( \frac{\partial p(\Ov{\vr}, \Ov{\vt}) }{\partial \vt} -
\frac{\partial p(\Ov{\vr}_\ep , \Ov{\vt}) }{\partial \vt} \right) \frac{\vte - \Ov{\vt} }{\ep} \Delta \Phi_\ep } \ \dt
= \chi(\ep, \eta).
\]

Consequently, we get
\bFormula{1---}
- \frac{1}{\ep} \int_0^\tau \intO{ \left[ \vre \Big( s(\vre,\vte)
- s(r_\ep, \Theta_\ep) \Big) \partial_t T_\ep + \vre \Big(
s(\vre,\vte) - s(r_\ep, \Theta_\ep) \Big) \vc U_\ep \cdot \Grad
T_\ep \right] } \ \dt
\eF
$$
- \frac{1}{\ep^2} \int_0^\tau \intO{ \frac{\vr_\ep - r_\ep}
{r_\ep}\partial_t p(r_\ep , \Theta_\ep )  } \ \dt -\frac 1{\ep^2}
\int_0^\tau
\intO{\Big(p(\vre,\vte)-p(\overline\vr_\ep,\overline\vt)\Big)\Delta\Phi_\ep}\ {\rm
d}t
$$
$$
= \int_0^\tau\intO{\Big(\delta T_\ep-\beta R_\ep\Big)\partial_t
T_\ep}{\rm d} t + \int_0^\tau\intO{R_\ep\partial_t\Big(\alpha
R_\ep+\beta T_\ep\Big)}\ {\rm d}t
$$
$$
-\left[ \int_0^\tau\intO{\Big(\delta\frac{\vte-\overline\vt}\ep
-\beta\frac{\vr_\ep-\overline\vr_\ep}\ep\Big)\partial_t T_\ep}\ {\rm d}t
+
\int_0^\tau\intO{\Big(\frac{\beta^2}{\beta^2+\alpha\delta}\frac{\vre-\overline\vr_\ep}\ep
-\frac{\beta\delta}{\beta^2+\alpha\delta}\frac{\vte-\overline\vt}\ep\Big)
\partial_t\Big(\alpha
R_\ep+\beta T_\ep\Big)}\ {\rm d}t \right]
$$
$$
- \int_0^\tau\intO{\Big(\delta\frac{\vte-\Theta_\ep}\ep
-\beta\frac{\vr_\ep-r_\ep}\ep\Big)\vc U_\ep\cdot\Grad T_\ep} \
{\rm d}t + \chi(\ep, \eta).
$$

In the next steps, we use the identities
\bFormula{ident}
(\beta^2+\alpha\delta) T=\beta(\alpha R+\beta T)+\alpha(\delta
T-\beta R),
\
(\beta^2+\alpha\delta) R=\delta(\alpha R+\beta T)-\beta(\delta
T-\beta R),
\eF
to compute,
\begin{equation}\label{*2}
\int_0^\tau\intO{\Big(\delta T_\ep-\beta R_\ep\Big)\partial_t
T_\ep}\ {\rm d} t + \int_0^\tau\intO{R_\ep\partial_t\Big(\alpha
R_\ep+\beta T_\ep\Big)}\ {\rm d}t
\end{equation}
$$
=
\int_{0}^\tau\int_{\Omega}\Big[\frac\beta{\beta^2+\alpha\delta}\Big(\delta
T_\ep-\beta R_\ep\Big)\partial_t\Big(\alpha R_\ep+\beta T_\ep\Big)
+ \frac\alpha{\beta^2+\alpha\delta}\Big(\delta T_\ep-\beta
R_\ep\Big)\partial_t\Big(\delta T_\ep-\beta R_\ep\Big)
$$
$$
+\frac\delta{\beta^2+\alpha\delta}\Big(\alpha R_\ep+\beta
T_\ep \Big)\partial_t\Big(\alpha R_\ep+\beta T_\ep \Big) -
\frac\beta{\beta^2+\alpha\delta}\Big(\delta T_\ep-\beta
R_\ep\Big)\partial_t\Big(\alpha R_\ep+\beta T_\ep\Big)
\Big]{\rm d}x\ {\rm d}t
$$
$$
=\frac12\frac\delta{\beta^2+\alpha\delta} \left[\int_{\Omega}|\alpha
R_\ep+\beta T_\ep|^2{\rm d} x\right]_0^\tau+
\frac12\frac\alpha{\beta^2+\alpha\delta} \left[\int_{\Omega}|\delta
T_\ep-\beta R_\ep|^2{\rm d} x\right]_0^\tau ,
$$
where we have used (\ref{r3}).

Similarly, we get
\begin{equation}\label{*3}
- \int_0^\tau\intO{\Big(\delta\frac{\vte-\overline\vt}\ep
-\beta\frac{\vr_\ep-\overline\vr}\ep\Big)\partial_t T_\ep}\ {\rm d}t
-
\int_0^\tau\intO{\Big(\frac{\beta^2}{\beta^2+\alpha\delta}\frac{\vre-\overline\vr}\ep
-\frac{\beta\delta}{\beta^2+\alpha\delta}\frac{\vte-\overline\vt}\ep\Big)
\partial_t\Big(\alpha
R_\ep+\beta T_\ep\Big)}\ {\rm d}t
\end{equation}
$$
=
-\frac{\alpha}{\beta^2+\alpha\delta}\int_0^\tau\intO{\Big(\delta\frac{\vte-\overline\vt}\ep-\beta\frac{\vre-\overline\vr}\ep\Big)
\partial_t\Big(\delta T_\ep-\beta R_\ep\Big)}\ {\rm d} t
$$

Finally, the last line on the right-hand side of (\ref{1---}) reads
\begin{equation}\label{*4}
-\int_0^\tau\intO{\Big(\delta\frac{\vte-\Theta_\ep}\ep
-\beta\frac{\vr_\ep-r_\ep}\ep\Big)\vc U_\ep\cdot\Grad T_\ep}\ {\rm
d}t
\end{equation}
$$=
-\frac\beta{\beta^2+\alpha\delta}\int_0^\tau\intO{\Big(\delta\frac{\vte-\Theta_\ep}\ep
-\beta\frac{\vr_\ep-r_\ep}\ep\Big)\vc U_\ep\cdot\Grad\Big( \alpha
R_\ep +\beta T_\ep\Big)}\ {\rm d}t
$$
$$
-
\frac\alpha{\beta^2+\alpha\delta}\int_0^\tau\intO{\Big(\delta\frac{\vte-\Theta_\ep}\ep
-\beta\frac{\vr_\ep-r_\ep}\ep\Big)\vc U_\ep\cdot\Grad\Big( \delta
T_\ep -\beta R_\ep\Big)}\ {\rm d}t
$$
$$
= -
\frac\alpha{\beta^2+\alpha\delta}\int_0^\tau\intO{\Big(\delta\frac{\vte-\Theta_\ep}\ep
-\beta\frac{\vr_\ep-r_\ep}\ep\Big)\vc U_\ep\cdot\Grad\Big( \delta
T_\ep -\beta R_\ep\Big)}\ {\rm d}t
 +\chi(\ep,\eta),
$$
where we have used  the dispersive estimates (\ref{a2}).

Summing up the previous integrals and using equation (\ref{r4a}) we may infer that
\begin{equation}\label{r12}
- \frac{1}{\ep} \int_0^\tau \intO{ \left[ \vre \Big( s(\vre,\vte)
- s(r_\ep, \Theta_\ep) \Big) \partial_t T_\ep + \vre \Big(
s(\vre,\vte) - s(r_\ep, \Theta_\ep) \Big) \vc U_\ep \cdot \Grad
T_\ep \right] } \ \dt
\end{equation}
$$
- \frac{1}{\ep^2} \int_0^\tau \intO{ \frac{\vr_\ep - r_\ep}
{r_\ep}\partial_t p(r_\ep , \Theta_\ep )  } \ \dt -\frac 1{\ep^2}
\int_0^\tau
\intO{\Big(p(\vre,\vte)-p(\overline\vr,\overline\vt)\Big)\Delta\Phi_\ep}\ {\rm
d}t
$$
$$
= \frac12\frac\delta{\beta^2+\alpha\delta} \left[\intO{|\alpha
R_\ep+\beta T_\ep|^2}\right]_0^\tau+
\frac12\frac\alpha{\beta^2+\alpha\delta} \left[\intO{|\delta
T_\ep-\beta R_\ep|^2}\right]_0^\tau
$$
$$
-\frac\beta{\beta^2+\alpha\delta}\int_0^\tau\intO{\Big(\delta\frac{\vte-\overline\vt}\ep
-\beta\frac{\vr_\ep-\overline\vr}\ep \Big)\vc v\cdot\Grad F_\ep}{\rm d}t +\chi(\ep,\eta)
$$

Finally, we use relation (\ref{Z1}) to obtain that
\[
-\frac\beta{\beta^2+\alpha\delta}\int_0^\tau\intO{\Big(\delta\frac{\vte-\overline\vt}\ep
-\beta\frac{\vr_\ep-\overline\vr}\ep \Big)\vc v\cdot\Grad F_\ep}{\rm d}t = \int_0^\tau \intO{ \frac{\vre - \Ov{\vr}_\ep}{\ep} \vc{v} \cdot \Grad F } \ \dt +
\chi(\ep, \eta),
\]
while { due to (\ref{r4a}) and (\ref{brum})}
\[
\frac12\frac\alpha{\beta^2+\alpha\delta} \left[\intO{|\delta
T_\ep-\beta R_\ep|^2}\right]_0^\tau = \frac{\beta}{\alpha} \int_0^\tau \intO{ T_\eta \vc{v} \cdot \Grad F } \ \dt + \chi(\ep, \eta).
\]

As $\theta$ and $T_{\eta}$ satisfy the \emph{same} transport equation and the acoustic system (\ref{r3}), (\ref{r4}) conserves the total energy,
we may use the previous estimates to rewrite (\ref{r11+}) in the final form:
\bFormula{rfinal}
\mathcal{E}_\ep \left( \vre, \vte, \vue \Big| r_\ep , \Theta_\ep , \vc{U}_\ep \right) (\tau) \leq \chi(\ep, \eta)  +
c \int_0^\tau  \mathcal{E}_\ep \left( \vre, \vte, \vue \Big| r_\ep , \Theta_\ep , \vc{U}_\ep \right) \ \dt,
\eF
which, performing the limit (i) for $\ep \to 0$, and then (ii) $\eta \to 0$, yields the conclusion of Theorem \ref{Tm1}.


\begin{thebibliography}{10}

\bibitem{AL2}
T.~Alazard.
\newblock Incompressible limit of the nonisentropic {E}uler equations with the
  solid wall boundary conditions.
\newblock {\em Adv. Differential Equations}, 10(1):19--44, 2005.

\bibitem{AL1}
T.~Alazard.
\newblock Low {M}ach number flows and combustion.
\newblock {\em SIAM J. Math. Anal.}, 38(4):1186--1213 (electronic), 2006.

\bibitem{AL}
T.~Alazard.
\newblock Low {M}ach number limit of the full {N}avier-{S}tokes equations.
\newblock {\em Arch. Rational Mech. Anal.}, {\bf 180}:1--73, 2006.

\bibitem{EF100}
E.~Feireisl.
\newblock Low {M}ach number limits of compressible rotating fluids.
\newblock {\em J. Math. Fluid Mechanics}, {\bf 14}:61--78, 2012.

\bibitem{FeNo6}
E.~Feireisl and A.~Novotn{\' y}.
\newblock {\em Singular limits in thermodynamics of viscous fluids}.
\newblock Birkh{\" a}user-Verlag, Basel, 2009.

\bibitem{FeiNov10}
E.~Feireisl and A.~Novotn{\' y}.
\newblock Weak-strong uniqueness property for the full
  {N}avier-{S}tokes-{F}ourier system.
\newblock {\em Arch. Rational Mech. Anal.}, {\bf 204}:683--706, 2012.

\bibitem{FeiNov12}
E.~Feireisl and A.~Novotn{\' y}.
\newblock Inviscid incompressible limits of the full
  {N}avier-{S}tokes-{F}ourier system.
\newblock {\em Commun. Math. Phys.}, {\bf 321}:605--628, 2013.

\bibitem{FeNoSun2}
E.~Feireisl, A.~Novotn{\' y}, and Y.~Sun.
\newblock { Dissipative solutions and the incompressible inviscid limits of the compressible magnetohydrodynamic system in unbounded domains.
\newblock {\em Disc. Cont. Dyn. Syst.} 34(1):121-143, 2014.}

\bibitem{FeiSch}
E.~Feireisl and M.E. Schonbek.
\newblock On the {O}berbeck-{B}oussinesq approximation on unbounded domains.
\newblock In {\em Abel Symposium Lecture Notes}. Springer Verlag, Berlin, 2011.

\bibitem{Isoz}
H.~Isozaki.
\newblock Singular limits for the compressible {E}uler equation in an exterior
  domain.
\newblock {\em J. Reine Angew. Math.}, 381:1--36, 1987.

\bibitem{JeJiNo}
D.~Jessl{\' e}, B.J. Jin, and A.~Novotn{\' y}.
\newblock {N}avier-{S}tokes-{F}ourier system on unbounded domains: weak
  solutions, relative entropies, weak-strong uniqueness. { {\em SIAM J. Math. Anal.}, 2013,}
\newblock to appear


\bibitem{Kato}
T.~Kato.
\newblock Remarks on the zero viscosity limit for nonstationary
  {N}avier--{S}tokes flows with boundary.
\newblock {\em In Seminar on PDE's, S.S. Chern (ed.), Springer, New York},
  1984.

\bibitem{KM1}
S.~Klainerman and A.~Majda.
\newblock Singular limits of quasilinear hyperbolic systems with large
  parameters and the incompressible limit of compressible fluids.
\newblock {\em Comm. Pure Appl. Math.}, {\bf 34}:481--524, 1981.

\bibitem{KL1}
R.~Klein.
\newblock Asymptotic analyses for atmospheric flows and the construction of
  asymptotically adaptive numerical methods.
\newblock {\em Z. Angw. Math. Mech.}, {\bf 80}:765--777, 2000.

\bibitem{Klein}
R.~Klein.
\newblock Scale-dependent models for atmospheric flows.
\newblock In {\em Annual review of fluid mechanics. {V}ol. 42}, Annu. Rev.
  Fluid. Mech., pages 249--274. Annual Reviews, Palo Alto, CA, 2010.

\bibitem{MAS5}
N.~Masmoudi.
\newblock Incompressible inviscid limit of the compressible {N}avier--{S}tokes
  system.
\newblock {\em Ann. Inst. H. Poincar{\' e}, Anal. non lin\' eaire}, {\bf
  18}:199--224, 2001.

\bibitem{MAS1}
N.~Masmoudi.
\newblock Examples of singular limits in hydrodynamics.
\newblock {\em In Handbook of Differential Equations, III, C. Dafermos, E.
  Feireisl Eds., Elsevier, Amsterdam}, 2006.

\bibitem{ZEY1}
R.~Kh. Zeytounian.
\newblock {J}oseph {B}oussinesq and his approximation: a contemporary view.
\newblock {\em C.R. Mecanique}, {\bf 331}:575--586, 2003.

\end{thebibliography}

\def\cprime{$'$} \def\ocirc#1{\ifmmode\setbox0=\hbox{$#1$}\dimen0=\ht0
  \advance\dimen0 by1pt\rlap{\hbox to\wd0{\hss\raise\dimen0
  \hbox{\hskip.2em$\scriptscriptstyle\circ$}\hss}}#1\else {\accent"17 #1}\fi}

\end{document}